\documentclass[a4paper]{jpconf}
\usepackage{graphicx}
\usepackage{amsmath,amsthm,amssymb,amsopn,amsfonts}
\usepackage{mathtools}
\usepackage{multirow}
\usepackage{xspace}
\usepackage{xcolor,colortbl}
\usepackage{booktabs}

\newcommand{\kron}{\otimes}
\newcommand{\eg}{e.g.\xspace}
\newcommand{\ie}{i.e.\xspace}
\renewcommand{\etal}{\textit{et al}.\xspace}

\newcommand{\R}{\mathbb{R}}
\newcommand{\p}[1]{\left( #1 \right)}
\newcommand{\mexp}[1]{\exp{\left( #1 \right)}}

\newcommand{\degree}{^{\circ}}

\newcommand{\TV}{TV\xspace}
\newcommand{\PMTV}{WF-TV\xspace}
\newcommand{\PVTV}{NLM-TV\xspace}
\newcommand{\PMFBP}{WF-FBP\xspace}
\newcommand{\PVFBP}{NLM-FBP\xspace}
\newcommand{\ang}{\AA\xspace}
\begin{document}

\bibliographystyle{iopart-num}

\title{Low-rank flat-field correction for artifact reduction in spectral computed tomography}

\author{Katrine O Bangsgaard$^1$, Genoveva Burca$^{2,3}$, Evelina Ametova$^{3,4}$, Martin S Andersen$^1$ and Jakob S J\o{}rgensen$^{1,3}$}

\address{
$^1$Department of Applied Mathematics and Computer Science, Technical University of Denmark, Richard Petersens Plads, Building 324, 2800 Kongens Lyngby, Denmark\\
$^2$ISIS Pulsed Neutron and Muon Source, STFC, UKRI, Rutherford Appleton Laboratory, Didcot OX11 0QX, United Kingdom\\
$^3$Department of Mathematics, The University of Manchester, Oxford Road, Alan Turing Building, Manchester M13 9PL, United Kingdom\\
$^4$Laboratory for Applications of Synchrotron Radiation, Karlsruhe Institute of Technology, Karlsruhe 76131, Germany
}

\ead{\{kaott, mskan, jakj\}@dtu.dk}

\begin{abstract}
Spectral computed tomography has received considerable interest in recent years since spectral measurements contain much richer information about the object of interest. In spectral computed tomography, we are interested in the energy channel-wise reconstructions of the object. However, such reconstructions suffer from low signal-to-noise ratio and share the challenges of conventional low-dose computed tomography such as ring artifacts. Ring artifacts arise from errors in the flat-field correction and can significantly degrade the quality of the reconstruction. We propose an extended flat-field model that exploits high correlation in the spectral flat-fields to reduce ring artifacts in the channel-wise reconstructions. The extended model relies on the assumption that the spectral flat-fields can be well-approximated by a low-rank matrix. Our proposed model works directly on the spectral flat-fields and can be combined with any existing reconstruction model, \eg, filtered back projection and iterative methods. The proposed model is validated on a neutron data set. The results show that our method successfully diminishes ring artifacts and improves the quality of the reconstructions. Moreover, the results indicate that our method is robust; it only needs a single spectral flat-field image, whereas existing methods need multiple spectral flat-field images to reach a similar level of ring reduction.
\end{abstract}

\section{Introduction}
Computed Tomography (CT) is a non-invasive imaging technique that allows us to obtain structural knowledge about the interior of objects from a set of projection images. Projection images are acquired by illuminating the object from different angles with radiation from a source \eg, an X-ray beam or beam of neutron radiation. The beam is attenuated as it travels through the object and the attenuated beam is measured by a detector placed opposite the source. {The attenuation is governed by absorption in X-ray CT and by scattering in the case of neutron CT. In both cases, the attenuation is material- and energy-specific, and if we measure the attenuation for multiple energies, also referred to as spectral CT, we can obtain a material decomposition of the object.}

To obtain the material decomposition, we need to compute the energy channel-wise reconstructions. Most reconstruction methods rely on the assumption that the detector response is known. In practice, however, the detector response is subject to various errors and must be estimated from measurements acquired without an object in the scanner, i.e., from flat-fields, also referred to as air scans \cite{Whiting2006}, white fields \cite{Sijbers2004} or open beams. {The flat-fields are noisy due to factors such as measurement noise, miscalibration, defective pixel elements with non-linear response, dust on scintillators}, and may introduce concentric rings in the reconstruction, also known as ring artifacts \cite{Goldman2007}. Ring artifacts are a great challenge for experimental CT set-ups with low-dose and/or short exposure time \cite{Aggrawal2018} and can significantly degrade the quality of the reconstruction. In spectral CT, we measure spectral flat-fields, \ie, flat-fields for each energy. However, the spectral measurements share the characteristics of low-dose CT since each energy channel has a low signal-to-noise ratio (SNR) and thus ring artifacts present a challenge in spectral CT \cite{Trapani2021, Zeng2016}.

\begin{figure}[ht!]
\centering
\includegraphics[width=0.7\textwidth]{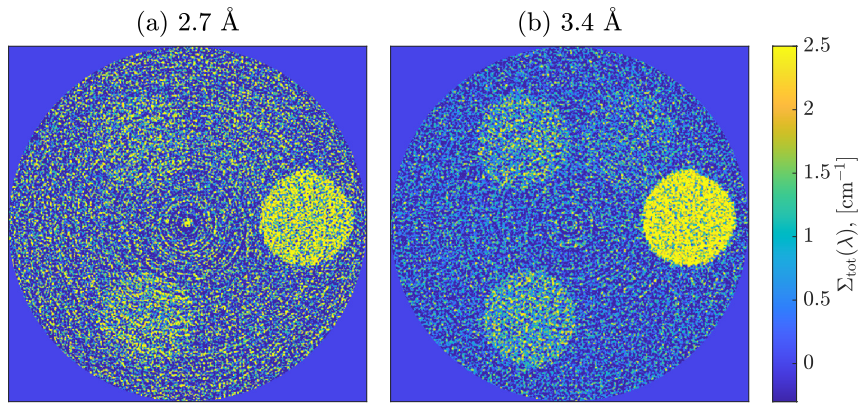}
\caption{FBP reconstructions of energy 2.7 \ang and 3.4 \ang. Ring artifacts are apparent in both reconstructions with varying severity.}
\label{fig:FBPrec}
\end{figure}

To illustrate the challenges of spectral CT, let us consider a neutron CT data set \cite{Neutrondata} which is described in detail in Section \ref{sec:Neutrondata}. Filtered back projection (FBP) reconstructions of the neutron data are shown for two energies in Figure \ref{fig:FBPrec} and the reconstructions reveal presence of ring artifacts. 

\subsection{Existing methods for ring reduction}
{Several reconstruction methods have been proposed to combat ring artifacts as part of the spectral reconstruction step. Wu \etal \cite{Dufan2015} propose a reconstruction method that exploits the similarity across spectral images by computing a polychromatic reconstruction (average across the spectral dimension) as a reference image combined with total variation (TV). Lv \etal \cite{Lv2020} and Fang \etal \cite{fang2020a} both propose deep learning approaches to suppress noise and remove ring artifacts for spectral CT. However, all methods rely on computationally expensive algorithms where the ring reduction is part of the reconstruction estimation process. 

Conventional preprocessing methods for ring reduction in monochromatic CT, \ie, single energy CT, can also be applied to the spectral CT data \cite{Mnch2009,Vo2018}. The main drawback here is that the ring reduction techniques must be applied channel-wise to the measured data.}

\subsection{Contribution}
\begin{figure}[ht!]
    \centering
    \includegraphics[width=0.75\textwidth]{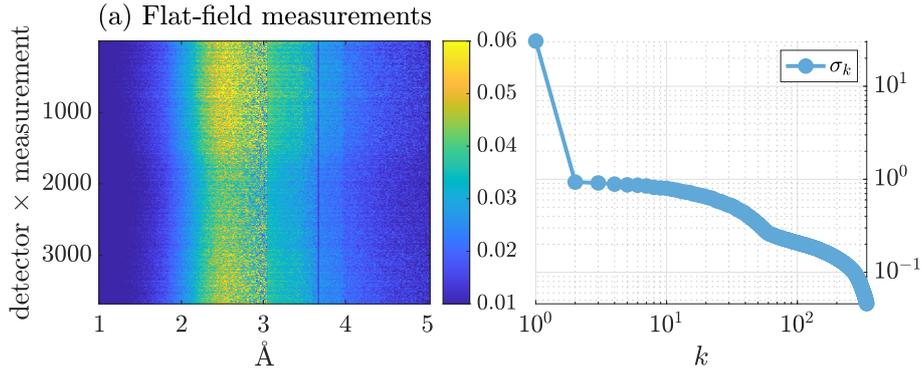}
    \caption{Visualization of the eight spectral flat-fields stacked vertically (a) and a loglog plot of the singular values (b). The singular values indicate that the spectral flat-fields are well-approximated by a rank-one matrix. }
    \label{fig:SVDintro}
\end{figure}

Figure \ref{fig:SVDintro} shows the eight measured spectral flat-fields stacked vertically and the corresponding singular values. A visual inspection suggests that the spectral flat-fields carry significant redundant information and that we can improve the SNR level in the spectral flat-fields by approximating the spectral flat-fields with a low-rank matrix. In particular, the singular values indicate that the spectral flat-fields can be well-approximated by a rank-one matrix due to the large jump in magnitude between the first and second singular values. {A similar idea where principal component analysis (PCA) is used to reduce ring artifacts in case of beam instability has been proposed by Hagemann \etal \cite{Hagemann2021} and Nieuwenhove \etal \cite{VanNieuwenhove2015}. The underlying assumptions in these studies are related, but the nature of the problems solved differs. }

Inspired by Figure \ref{fig:SVDintro}, we propose an extended flat-field model that exploits high correlation across channels in the spectral flat-fields to reduce ring artifacts in the reconstructions. The extended model relies on the assumption that the spectral flat-fields can be well-approximated by a low-rank matrix. Our method does not depend on a specific reconstruction method since it works directly on the spectral flat-fields. Hence, our method can be combined with all existing reconstruction models such as the conventional FBP or more advanced spectral methods e.g., \cite{Dufan2015,Zhang2016,Gao2011,Ametova2021}. {Moreover, our method does not need to be applied channel-wise in the sense that the low-rank spectral flat-field simply replaces the measured spectral flat-field in the reconstruction step.}

\subsection{Outline}
Section 2 introduces the spectral CT model, existing methods for ring reduction and the proposed methodology. In Section 3, the experimental set-up for the neutron data set is described and numerical experiments are conducted in Section 4. Section 5 discusses the results and Section 6  concludes the paper.

\subsection{Notation}
The set $\R^n$ is the $n$-dimensional Euclidean space, $\R_+^n$ is the non-negative orthant and $\R^{m \times n}$ denotes the set of real-valued $m\times n$ matrices. The vector $\textbf{1}_n\in \R^n$ is a vector of ones, $I_{n\times n}\in \R^{n\times n}$ denotes the identity matrix and the transpose of $A$ is denoted $A^T$. The function $\exp$ with a vector or matrix as argument are to be interpreted element-wise. If $x$ is a vector, then $\text{diag}(x)$ is the diagonal matrix with the elements of $x$ on the diagonal. The 2-norm of a vector is denoted $\| \cdot \|_2$ and the Frobenius norm is denoted $\|\cdot \|_{\text{F}}$.

\section{Methods}
Consider a spectral data set with $m$ energy channels and let $E_k$ denote the $k$th energy channel. The incident intensity of a beam with energy $E_k$ on a detector elemenet is prescribed by the Beer-Lambert law \cite{herman2009a}, 
\begin{align}\label{eq:lambertbeer}
I(E_k) = I_0(E_k) \mexp{-\int_{\ell} \mu(\textbf{x},E_k) \text{d}\textbf{x}} , 
\end{align}
where $I(E_k)$ and $I_0(E_k)$ are the energy-dependent intensity incident on the detector element and in the object, respectively. Further, $\ell$ is the line segment between the source and the detector, and $\mu\colon \R^d \times \R \to\R_+$ is the energy-dependent spatial attenuation function. 

Let $Y_k\in \R^{rp}$ denote the measurements for the $k$'th energy with $r$ detector elements and $p$ projection images and discretize the domain into $n$ pixels, then by appropriate discretization (\eg, see \cite{Bangsgaard2021}), we can describe the measurements by,
\begin{align}\label{eq:convenmodel}
    Y_k = \text{diag}\p{\textbf{1}_p \kron Z_k} \mexp{-A X_k},
\end{align} 
where $\kron$ is the Kronecker product, $Z_k\in \R^{r}$ is energy-dependent intensity incident on the detector, $X_k\in \R^n$ are the unknown attenuation coefficients and $A\in \R^{rp\times n}$ is the system matrix describing the traveled distance of the beams. In practice, the energy-dependent intensity incident on the detector is estimated by measuring the detector response. We assume that $s$ spectral flat-fields are measured for each energy. The estimate for $Z\in \R^{r\times m}$ is then conventionally computed by the mean of the $s$ measured spectral flat-fields, \ie,
\begin{align}\label{eq:convF}
\hat{Z}  = \frac{1}{s}\sum_{j=1}^s F_j= \frac{1}{s}\p{\textbf{1}_s^T\kron I_{r\times r}} F,
\end{align}
where $F = [F_1^T, F_2^T, \hdots , F_s^T]^T$ and $F_j\in \R^{r\times m}$ is the $j$'th spectral flat-field. However, if $F$ is noisy, then the estimate $\hat{Z}_k$ can give rise to ring artifacts in the reconstruction. 

\subsection{Low-rank approximation}
 Essentially, each of the spectral flat-fields carries information about the detector response for all energies and the aim is to exploit the high correlation in the spectral dimension motivated by the observations in Figure \ref{fig:SVDintro}. 
 
\begin{figure}[ht!]
    \centering
    \includegraphics[width=1\textwidth]{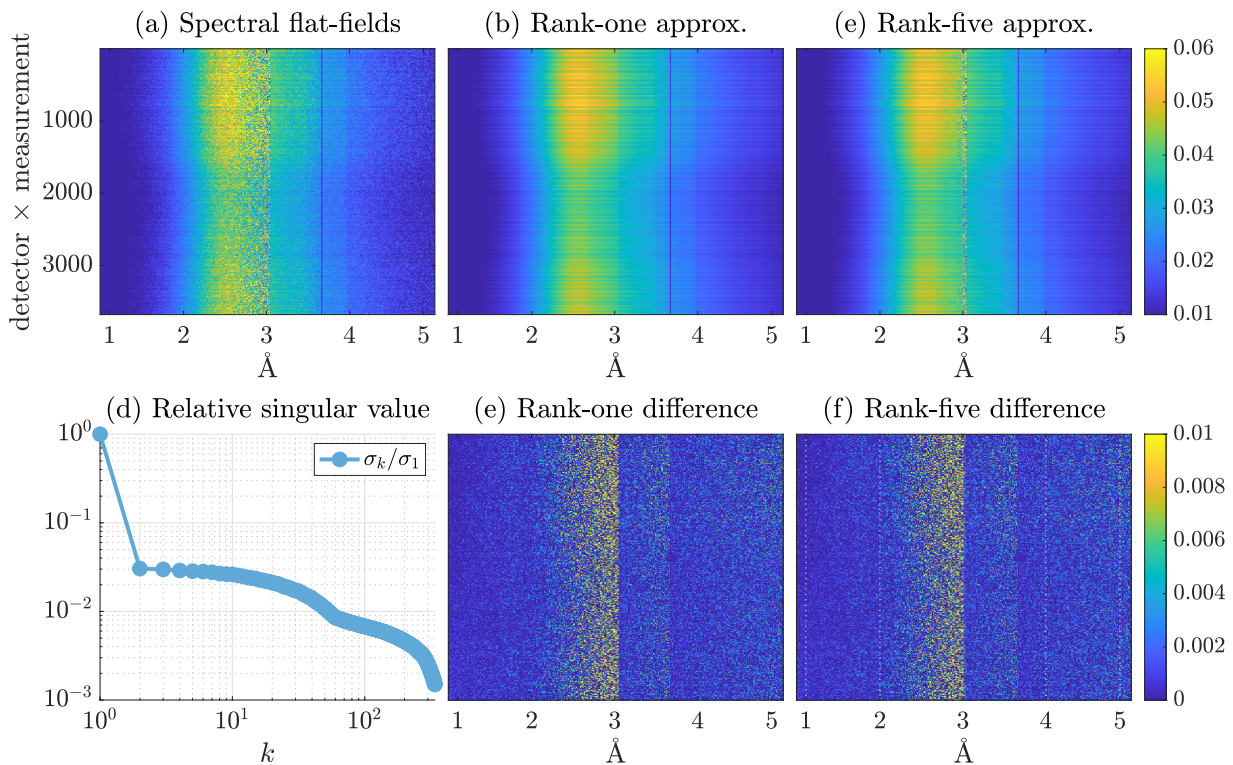}
    \caption{Rank-one and five approximations of the spectral flat-field are shown in (b) and (c). The corresponding difference images between the spectral flat-fields and low-rank matrices are shown in (e) and (f) for rank-one and rank-five, respectively.}
    \label{fig:SVDrank15}
\end{figure}

 The nearest low-rank matrix of $F$ in the spectral norm can be computed by means of a singular value decomposition (SVD) \cite{Loan2000}. An SVD of $F$ is a decomposition 
 \begin{align}
F = U \Sigma V^T,
\end{align}
where the columns of $U\in \R^{rs\times rs}$ and $V\in\R^{m\times m}$ are orthogonal matrices, and $\Sigma \in \R^{rs\times m}$ is a matrix with the singular values $\sigma_1  \ge \sigma_2  \ge \cdots  \ge 0$ on the main diagonal and zeroes elsewhere. The best rank-$l$ approximation of $F$ in the spectral norm is then given by
\begin{align}
F^l=  \sum_{i=1}^l\sigma_i U_i V_i^T,
\end{align}
where $U_i$ and $V_i$ are the $i$th columns of $U$ and $V$, respectively.
The relative approximation error is given by,
 \begin{align}\label{eq:approxerror}
\frac{\|F^l-F \|_{\text{F}}}{\|F\|_{\text{F}}}= \frac{\sigma_{l+1}}{\sigma_1},
 \end{align}
which follows from the Eckart--Young--Mirsky theorem.
The rank-one and rank-five approximations of the spectral flat-fields for the neutron data are shown in Figure \ref{fig:SVDrank15}. The approximation error for the rank-one and rank-five matrices are 0.030 and 0.028, respectively. Hence, we only reduce the approximation error with 0.002 by including four extra singular vectors. Considering Figure \ref{fig:SVDrank15}, we see that the rank-one approximation has mitigated a substantial amount of noise whereas in the rank-five approximation, some of the noisy tendencies in the spectral flat-fields start to reappear. This observation is confirmed by considering the difference images between the spectral flat-fields and the low-rank matrices in Figure \ref{fig:SVDrank15}e and \ref{fig:SVDrank15}f. Hence, for the numerical experiments, we will confine ourselves to considering the rank-one approximation of the spectral flat-fields, \ie, we compute the estimate for $\hat{Z}$ by replacing $F$ by $F^1$ in \eqref{eq:convF}, \ie,
\begin{align}\label{eq:Fnew}
   \hat{Z} = \frac{1}{s}\p{\textbf{1}_s^T\kron I_{r\times r}} F^1.
\end{align}

At this point, we emphasize that our method is computationally cheap as it scales with $\mathcal{O}(\min(rs,m)^2\max(rs,m))$. Moreover, it can be combined with any reconstruction method, e.g., FBP, iterative methods, statistical models, etc. since it is only applied to the spectral flat-fields. 

\subsection{Existing ring reduction method}\label{sec:Munch} 
{Our method works solely on the spectral flat-fields and thus separates the ring reduction from the reconstruction step. Hence, we compare our method to two existing ring reduction techniques for monochromatic CT.} We compare our methodology with the preprocessing method proposed by Münch \etal \cite{Mnch2009} which combines wavelet and Fourier filtering to mitigate ring artifacts in the reconstruction. The method is computationally inexpensive and does not increase the overall computational cost significantly. The method depends on three parameters; we use a damping factor of 0.9 and the Daubechies 5 wavelet with a three-level decomposition for all numerical experiments, see \cite{Mnch2009} for details. The second preprocessing method is proposed by Vo \etal \cite{Vo2018} and uses a combination of sorting and smoothing (non-local means) in attempt to smooth the data and thereby reduce ring artifacts. The method depends on a parameter related to the smoothing filter, we choose a parameter value of $31$, see \cite{Vo2018} for further details. We denote the methods from Münch \etal \cite{Mnch2009} and Vo \etal \cite{Vo2018} by FW (Wavelet Fourier) and NLM (non-local means), respectively.


\section{Neutron data}\label{sec:Neutrondata}
\begin{figure}[ht]\centering
\includegraphics[width=0.25\textwidth]{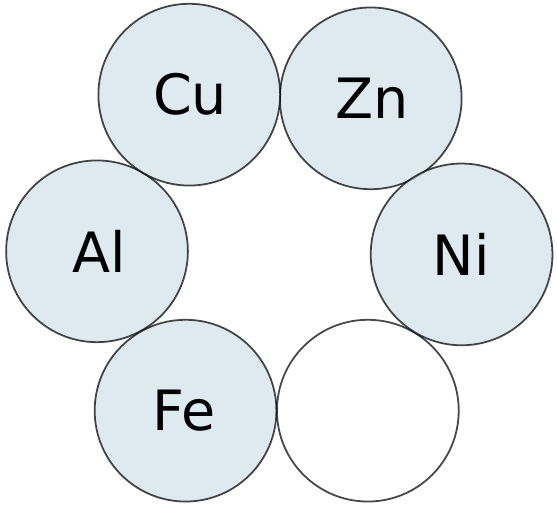}\hspace{0.05\textwidth}%
\begin{minipage}[b]{0.45\textwidth}\caption{\label{fig:sketchneutrondata}Sketch of the object for the experimental neutron data \cite{Neutrondata}. The experimental set up consists of six cylinders where five of the cylinders have been filled with aluminum (Al), iron (Fe), copper (Cu), nickel (Ni) and zinc (Zn) powders. The sixth cylinder is empty.}
\end{minipage}
\end{figure}

We validate the proposed methodology on a neutron CT data set \cite{Neutrondata}. The neutron data were acquired at the imaging and materials science and engineering (IMAT) beamline operating at the ISIS spallation neutron source (Rutherford Appleton Laboratory, UK). Figure \ref{fig:sketchneutrondata} shows a sketch of the object of interest. The object consists of six cylinders whereof five are filled with metal powders, \ie, aluminum (Al), iron (Fe), copper (Cu), nickel (Ni) and zinc (Zn) powders.

\subsection{Data acquisition}

 {Ametova \etal \cite{Ametova2021} describe the data acquisition and pre-treatment of data in detail, \ie, beam instabilities, overlap correction and spectral averaging. We consider the reconstruction-ready measurements and confine ourselves to reporting only the essential details of the experimental set-up.} The data set contains $m=339$ spectral projections acquired at $p=120$ equidistant angles distributed from $0\degree$ to $180\degree$ with $1.5\degree$ angular increments. Eight spectral flat-fields were acquired, four prior to the scan and four after, i.e., $s = 8$. Each projection consists of $460 \times 460 $ pixels with a  pixel size of  $0.055$ mm resulting in a view of approximately $25 \times 25$ mm$^2$ and $r = 460$ detector elements and $n = 460^2$ pixels for the experiments. We choose the $127$th vertical detector row for experiments, \ie, we consider a two-dimensional set-up. The experimental set-up is summarized in Table \ref{tab:Dimensions}. The attenuation coefficient for the neutron experiment is denoted $\Sigma_{\rm tot}(\lambda)$ and has unit cm$^{-1}$, see \cite{Ametova2021} for further details. 

\begin{table}[ht!]
\caption{\label{tab:Dimensions}Experimental set-up for the spectral neutron data used for the numerical experiments.}
\begin{center}
\begin{tabular}{*{5}c}
\br
pixels ($n$)&energy channels ($m$) &detectors ($r$) & projections ($p$) &  spectral flat-fields ($s$)\cr\mr
 $460^2$ & $339$ & $460$ & $120$ & $8$ \cr
\br
\end{tabular}
\end{center}
\end{table}

\section{Numerical experiments}
We consider two reconstruction models: FBP and a weighted least squares (WLS) reconstruction model combined with TV regularization. We compare our method to the conventional flat-field correction and the existing ring reduction techniques described in Subsection \ref{sec:Munch}. Table \ref{tab:Methods} provides an overview of the reconstruction models and ring reduction techniques.

\begin{table}[ht!]
    \centering
     \caption{\label{tab:Methods}Abbreviations used for the reconstruction models and ring reduction techniques.}
    \begin{tabular}{l l l l}\br
    \multirow{ 2}{*}{Ring reduction technique} & & \multicolumn{2}{c}{Reconstruction Model}    \cr
    & & FBP & WLS with TV \cr
         \cmidrule{1-1} \cmidrule{3-4}
    Conventional && FBP & TV  \cr
    Preprocessing Münch  \etal \cite{Mnch2009} &   & \PMFBP & \PMTV \cr
    Preprocessing Vo \etal \cite{Vo2018}   & & \PVFBP & \PVTV \cr
    Low-rank spectral flat-fields& & LR-FBP & LR-\TV \cr\br
    \end{tabular}
\end{table}

We use AIR TOOLS II \cite{Hansen2017} to generate the parallel-beam geometry of the experimental set-up. For the FBP reconstructions, we use the \texttt{fbp} function from AIR TOOLS II with the \texttt{Hann} filter to reduce the noise in the computed FBP reconstructions. For the TV reconstructions, we used the implementation of WLS with TV from \cite{Aggrawal2018} in MATLAB. All TV reconstructions have a maximum number of iterations of 1000 and a regularization parameter of $0.005$. The regularization parameter was found by visual inspection of reconstructions for varying values of the regularization parameter. { Note that we pick the same regularization parameter for all energies which will result in some reconstructions being a bit over-regularized whereas other reconstructions might be slightly under-regularized. The reason is that the SNR changes significantly as a function of energy and thus there is not a single regularization parameter that fits all energy channels. However, the purpose of the methods is ring reduction and thus we limit the experiments to considering the same regularization parameter for all reconstructions. }
\subsection{Error measures}
The quality of the computed reconstructions is assessed by visual inspection combined with contrast-to-noise ratio (CNR). The CNR metric is used for evaluating the image contrast and noise properties for a selected region of interest (ROI) \cite{warr2021a}. We use the method proposed by Bian \etal \cite{bian2010a} where a ROI with low-contrast structure is compared to a background ROI while taking the standard deviations of both the signal and background ROIs into account. The ROIs are clearly marked on the figures when applicable.

\subsection{Experiment: Different energies}
There are too many energies to visualize all reconstructions, and the SNR varies significantly between the energies \cite{Ametova2021}. Thus, we select three energies based on the relative difference (RD) between the computed FBP and LR-FBP solution to ensure a representative visualization of energy channels. We define RD by the measure,
\begin{align}\label{eq:relativediff}
   \text{RD}(k) =  \frac{\| X_k^{\text{FBP}}- X_k^{\text{LR-FBP}}\|_2}{\| X_k^{\text{LR-FBP}}\|_2}.
\end{align}
where $X^{\text{FBP}}_k$ and $X_k^{\text{LR-FBP}}$ denote the FBP and LR-FBP reconstructions for energy $k$, respectively. We select three energies corresponding to the minimum, median and maximum RD, corresponding to energies $3.1$ \ang, $4.2$ \ang and $2.9$ \ang, respectively.

\begin{table}[ht!]
\caption{\label{tab:Experiment1}CNR for experiment depicted in Figures \ref{fig:Experiment1FBP} and \ref{fig:Experiment1TV}. The ROIs are marked in Figures  \ref{fig:Experiment1FBP} and \ref{fig:Experiment1TV}. }
\begin{center}
\begin{tabular}{l*{10}l}
\br
\multirow{2}{*}{Energy} & &\multicolumn{4}{c}{FBP} && \multicolumn{4}{c}{\TV} \cr
 &  &\multicolumn{1}{c}{Conv.} & \multicolumn{1}{c}{WF} &\multicolumn{1}{c}{NLM} &  \multicolumn{1}{c}{LR} & &\multicolumn{1}{c}{Conv.} & \multicolumn{1}{c}{WF} &\multicolumn{1}{c}{NLM} &  \multicolumn{1}{c}{LR} \cr
\cmidrule{1-1}\cmidrule{3-6}\cmidrule{8-11}
3.1 \ang && 2.15&  2.31&  2.17&  2.31 &&    1.02$\cdot 10^{3}$  &   1.03$\cdot 10^{3}$  &  1.01$\cdot 10^{3}$  &  1.06$\cdot 10^{3}$ \cr
4.2 \ang && 0.28&  0.30&  0.29&  0.30 &&    0.16$\cdot 10^{2}$  &   0.16$\cdot 10^{3}$  &  0.21$\cdot 10^{3}$  &  0.28$\cdot 10^{3}$  \cr
2.9 \ang && 0.20&  0.24&  0.24&  0.25 &&    0.03$\cdot 10^{1}$  &   0.04$\cdot 10^{3}$  &  0.07$\cdot 10^{3}$  &  0.07$\cdot 10^{3}$  \cr
\br
\end{tabular}
\end{center}
\end{table}
    
\begin{figure}[hb!]
    \centering
    \includegraphics[width=1\textwidth]{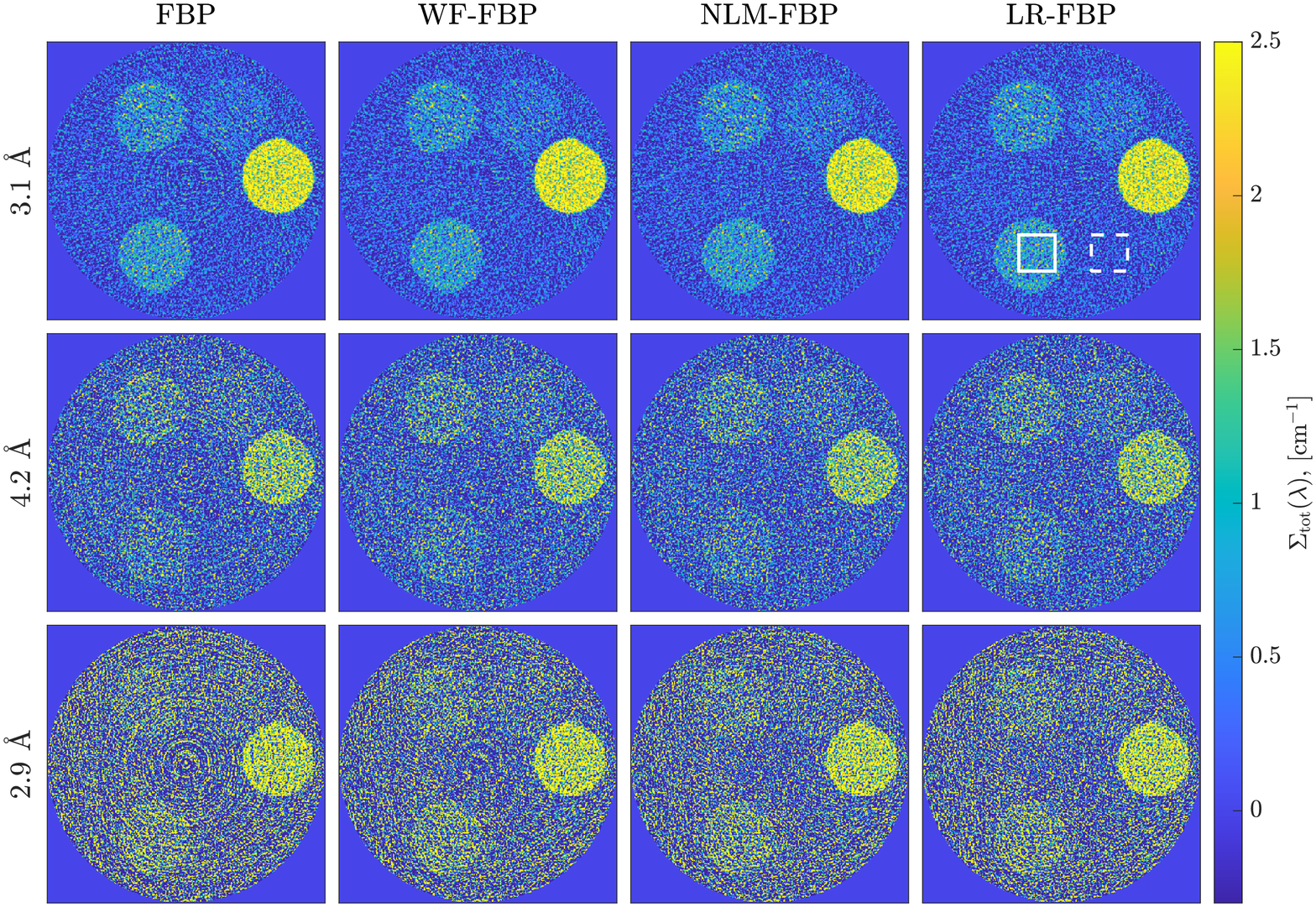}
    \caption{FBP (first column), \PMFBP (second column), \PVFBP (third column) and LR-FBP (fourth column) reconstructions for three energies chosen by the RD measure. The white squares with full and dashed lines mark the structure and background ROIs for CNR, respectively.}
    \label{fig:Experiment1FBP}
\end{figure}

\begin{figure}[ht!]
    \centering
    \includegraphics[width=1\textwidth]{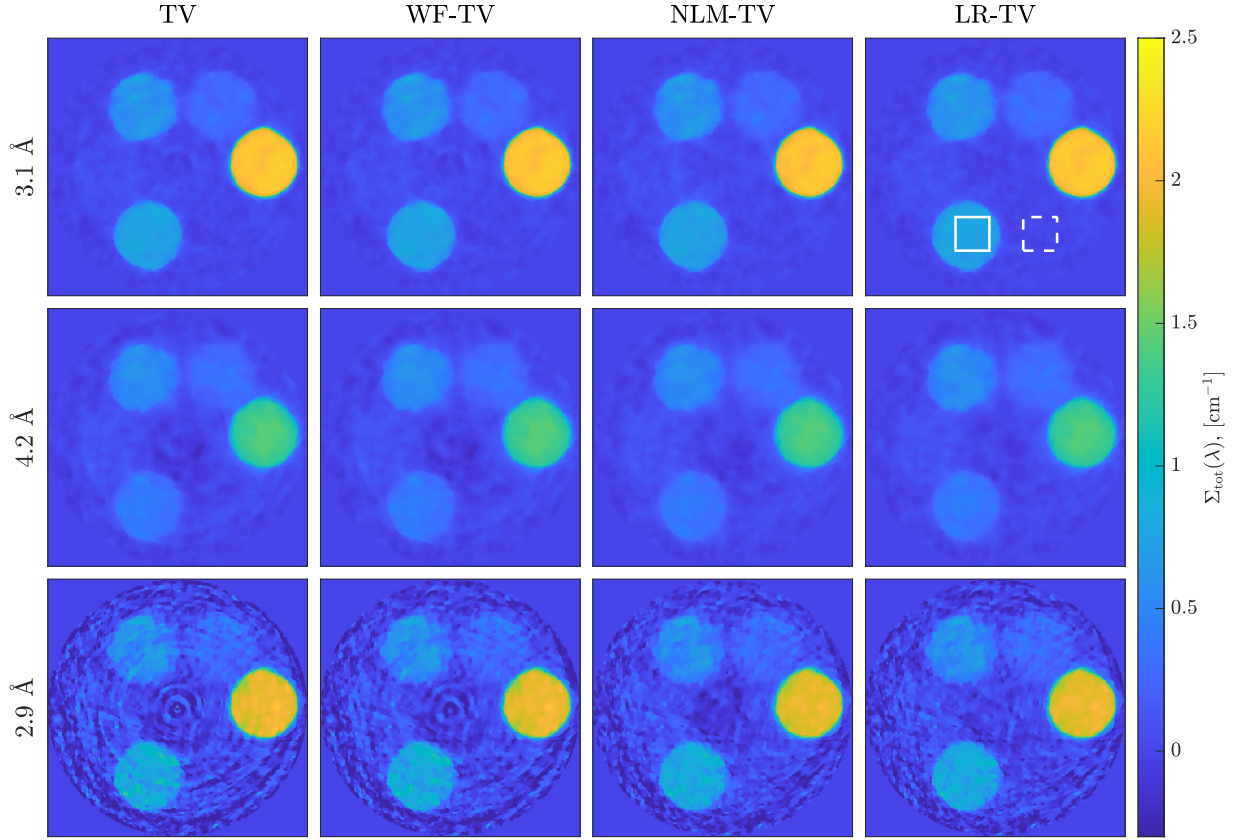}
    \caption{\TV (first column), \PMTV (second column), \PVTV (third column) and LR-TV (fourth column) reconstructions for three energies chosen by the RD measure. The white squares with full and dashed lines mark the structure and background ROIs for CNR, respectively.}
    \label{fig:Experiment1TV}
\end{figure}

The FBP reconstructions in Figure \ref{fig:Experiment1FBP} (first column) are all affected by ring artifacts. However, the severity of the ring artifacts increases from top to bottom, \ie, when the RD increases. The \PMFBP reconstructions (second column) seem to mitigate the ring artifacts at first glance, however, when carefully comparing the \PMFBP with \PVFBP and LR-FBP, wave-like artifacts can be seen in the \PMFBP reconstructions. Hence, the ring artifacts have been reduced but not eliminated. For the \PVFBP and LR-FBP reconstructions (third and fourth column), there are no visible artifacts remaining. 

Considering the CNR listed in Table \ref{tab:Experiment1}, we see that the CNR decreases from top to bottom for all four methods. Thus, when the SNR is low, our method differs most from the FBP, which can be explained by the fact that the ring artifacts are more dominating when the SNR is low. The CNR for all four methods is quite close in Figure \ref{fig:Experiment1FBP}. This might be explained by the fact that the dominating noise contribution comes from the measurement noise and not the ring artifacts, since all FBP reconstructions suffer from a very low SNR. Note that the high noise level might conceal remaining ring artifacts.

Inspecting the \TV reconstructions in Figure \ref{fig:Experiment1TV}, we see that TV regularization reduces the noise level significantly, and consequently, the ring artifacts appear more severe for the \TV reconstructions. In addition, we also note that the vague wave structures in the \PMFBP reconstructions are even clearer in the \PMTV reconstructions compared to the \PMFBP reconstructions. By carefully inspecting the \PVTV reconstruction for $2.9$ \ang, one can see that the preprocessing has introduced a dark spot with negative values in the center. The other \PVTV reconstructions show no introduced artifacts and closely resemble the LR-TV reconstructions (fourth column). The LR-\TV reconstructions reveal no ring structure even though the noise level is reduced. 

The TV reconstructions have significantly improved the overall CNR compared to the FBP reconstructions, which can be seen in Table \ref{tab:Experiment1}. For the TV experiment in Figure \ref{fig:Experiment1TV}, we also see an increase in CNR when applying a ring reduction method, especially for \PVTV and LR-TV. The LR-TV reconstruction achieves the highest CNR for two out of the three energies depicted in Figure \ref{fig:Experiment1TV}.

\begin{figure}[ht!]
    \centering
    \includegraphics[width=1\textwidth]{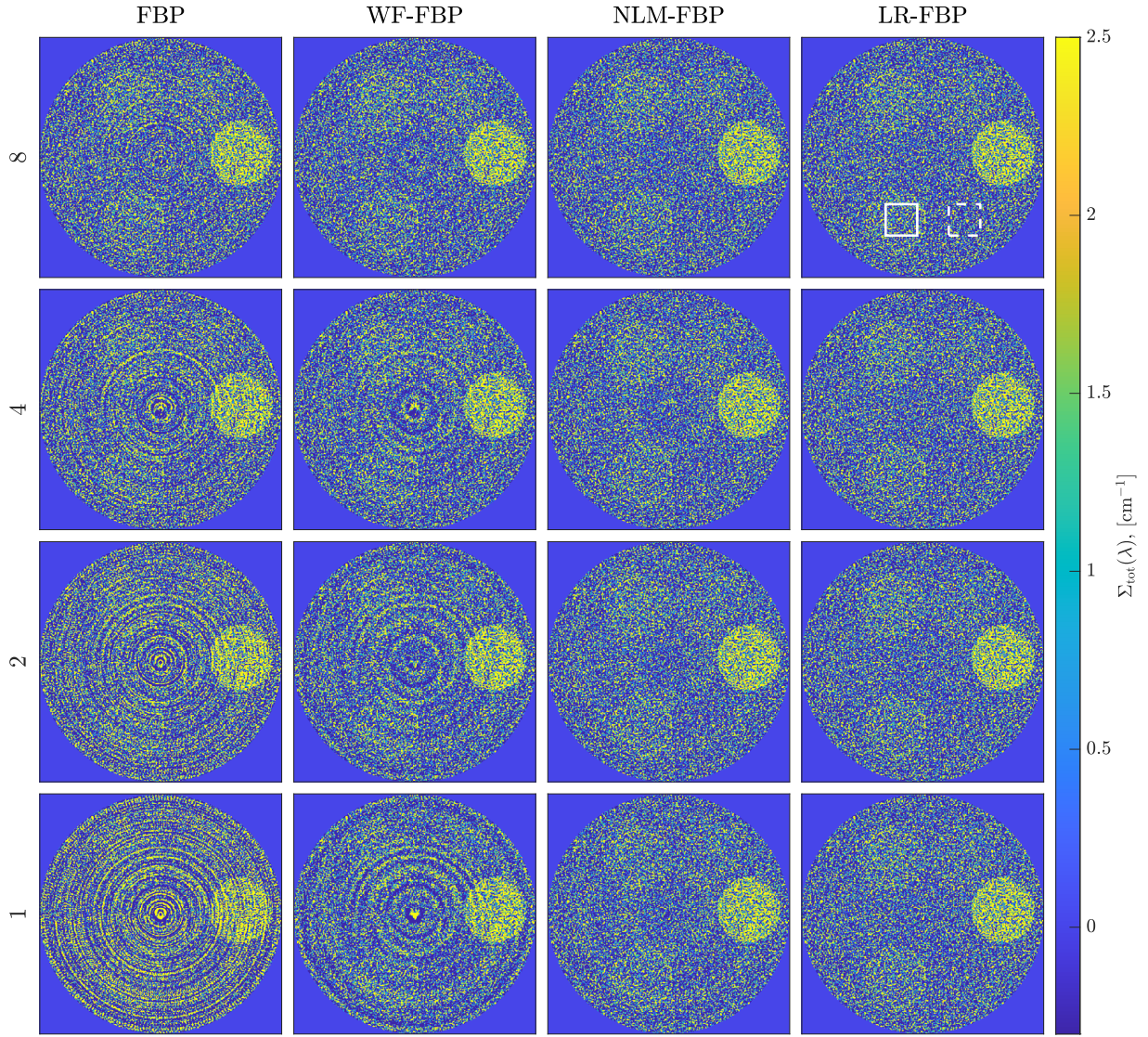}
    \caption{FBP, \PMFBP, \PVFBP and LR-FBP reconstructions of energy 2.1 \ang with eight, four, two and one flat-field, respectively. The white squares with full and dashed lines mark the structure and background ROIs, respectively. }
    \label{fig:Experiment2FBP}
\end{figure}

\begin{figure}[ht!]
    \centering
    \includegraphics[width=1\textwidth]{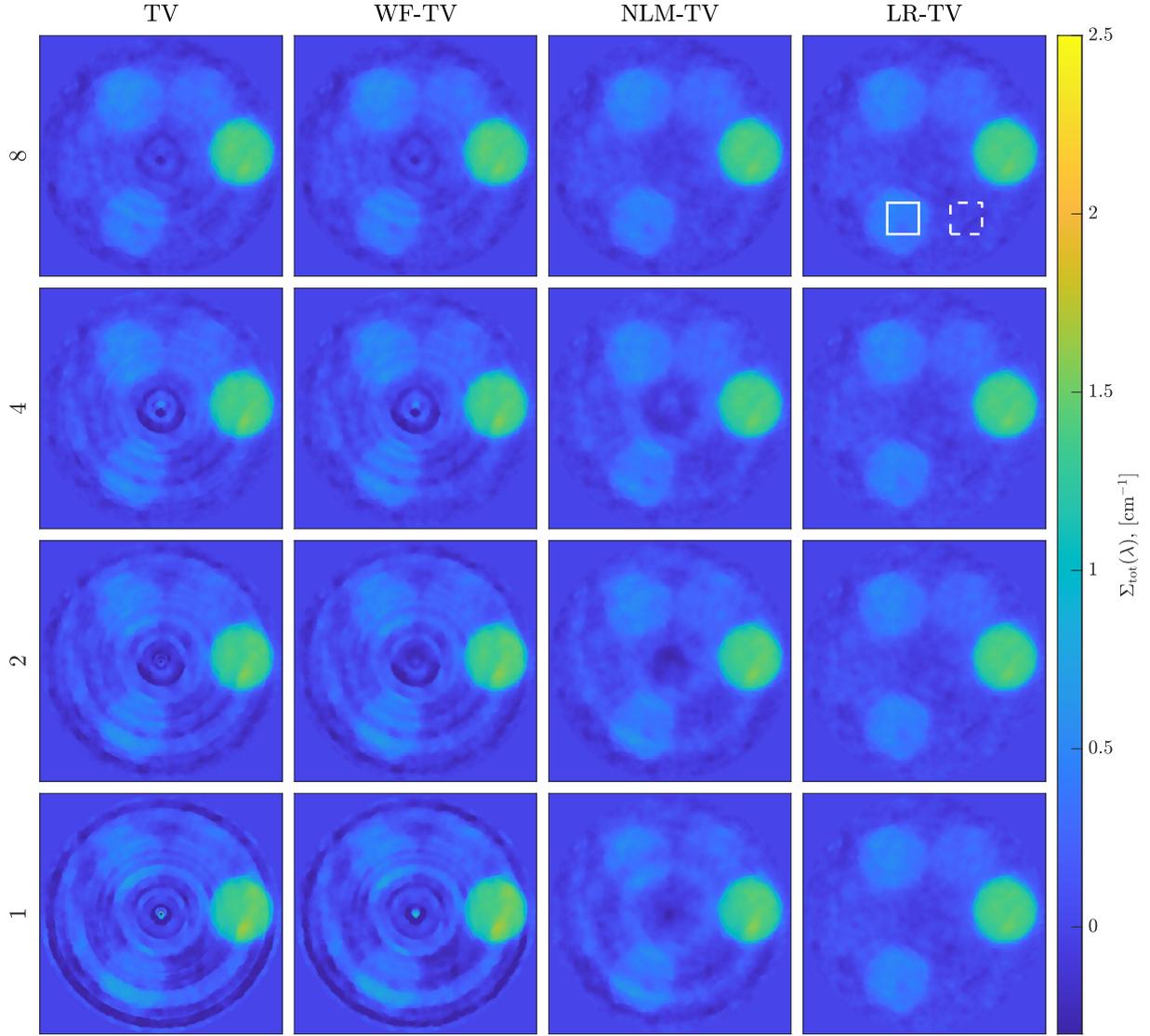}
    \caption{\TV, \PMTV, \PVTV and LR-\TV reconstructions of energy 2.1 \ang with eight, four, two and one flat-field, respectively. The white squares with full and dashed lines mark the structure and background ROIs, respectively. }
    \label{fig:Experiment2TV}
\end{figure}

\subsection{Experiment: Effect of number of spectral flat-fields}
We now perform an experiment with varying number of spectral flat-fields to validate the robustness of our proposed method. We consider energy 2.1 \ang for the experiment and use eight, four, two and one flat-field for the flat-field correction, respectively. The FBP, \PMFBP, \PVFBP and LR-FBP reconstructions are depicted in Figure \ref{fig:Experiment2FBP} and the \TV reconstructions for the same experiment are shown in Figure \ref{fig:Experiment2TV}. The CNR measures for the reconstructions are listed in Table \ref{tab:Experiment2}.

\begin{table}[ht!]
\caption{\label{tab:Experiment2}CNR for experiment depicted in Figure \ref{fig:Experiment2FBP} and Figure \ref{fig:Experiment2TV} for energy 2.1 \ang.}
\lineup
\begin{center}
\begin{tabular}{c*{11}c}
\br
\multirow{2}{*}{Flat-fields} & &\multicolumn{4}{c}{FBP} && \multicolumn{4}{c}{\TV} \cr
 &  &\multicolumn{1}{c}{Conv.} & \multicolumn{1}{c}{WF} &\multicolumn{1}{c}{NLM} &  \multicolumn{1}{c}{LR} & &\multicolumn{1}{c}{Conv.} & \multicolumn{1}{c}{WF} &\multicolumn{1}{c}{NLM} &  \multicolumn{1}{c}{LR} \cr
\cmidrule{1-1}\cmidrule{3-6}\cmidrule{8-11}
8 &&0.15 & 0.17 & 0.16 &  0.16 && 0.95$\times 10^{2}$ & 0.95$\times 10^{2}$ & 1.62$\times 10^{2}$ &  1.70$\times 10^{2}$\cr 
4 &&0.13 & 0.15 & 0.14 &  0.16 && 0.55$\times 10^{2}$ & 0.56$\times 10^{2}$ & 1.39$\times 10^{2}$ &  1.79$\times 10^{2}$\cr 
2 &&0.11 & 0.14 & 0.13 &  0.16 && 0.36$\times 10^{2}$ & 0.39$\times 10^{2}$ &  0.87$\times 10^{2}$ &  1.78$\times 10^{2}$\cr
1 &&0.06 & 0.12 & 0.14 &  0.16 && 0.14$\times 10^{2}$ & 0.15$\times 10^{2}$ &  0.52$\times 10^{2}$ &  1.80$\times 10^{2}$\cr  
\br
\end{tabular}
\end{center}
\end{table}

If we start by considering the first column of Figure \ref{fig:Experiment2FBP}, we see that there are visible ring artifacts in all FBP reconstructions, but the severity of the ring artifacts increases as the number of flat-fields used for flat-field correction decreases (\ie, from top to bottom). The \PMFBP reconstructions are shown in the second column of the figure. Vague ring artifacts can be seen in the \PMFBP reconstruction using all eight flat-fields. However, the severity of these ring artifacts increases when the number of flat-fields decreases, just as for the FBP reconstructions. The \PVFBP reconstructions in third column show less wave-like artifacts than \PVFBP, however, then using four or fewer flat-fields, the method starts to struggle, and artifacts start to arise in the reconstructions. The LR reconstructions in the fourth column show no sign of artifacts, not even in the case where a single flat-field is used and thus our method seems to be robust. The robustness of LR-FBP can also be seen in the CNR reported in Table \ref{tab:Experiment2}. Here the CNR is almost constant for the four reconstructions with LR, whereas the CNR increases with the number of flat-fields for both FBP, \PMFBP and \PVFBP. The experiment was repeated using TV for the reconstruction model and the findings support those found using FBP. The SNR is significantly increased by using TV regularization and the remaining artifacts for the preprocessing methods are even more pronounced in Figure \ref{fig:Experiment2TV}, which is also supported by the CNR measure listed in Table \ref{tab:Experiment2}.

\section{Discussion}
An explanation for the robustness of our method lies with the fact that the spectral dimension carries lots of redundant information. For the experiment, we have eight flat-fields but 339 energy channels. Thus, using a single spectral flat-field corresponds to 339 measurements for our method, whereas eight flat-fields with the conventional method only corresponds to eight flat-fields. {Thus, if the assumption of low rank holds, then more energy channels will give more redundancy in the spectral flat-fields and thereby our proposed method will gain further advantage over the existing methods.}

We investigated spectral plots for material decomposition, \ie, plot of a single pixel across the spectral dimension for the distinct materials. However, preliminary results showed no visual effect of applying the preprocessing methods. 

We could consider the spectral flat-fields as a three-dimensional tensor with dimensions $s \times r\times m$. While tensor decomposition methods exist, such as the Tucker decomposition and parallel factors decomposition (PARAFAC) \cite{Filipovi2013, Harshman1994,Grasedyck2013}, preliminary experiments showed no significant difference between the low-rank tensor, and the low-rank matrix obtained by treating the spectral flat-fields as a matrix as in \eqref{eq:convF}. 

 A possible improvement of the model could be the inclusion of the low-rank matrix as parameters of the reconstruction model and jointly estimating the reconstruction and the spectral detector response, \eg, see \cite{Bangsgaard2021}. A natural initial guess for the low-rank matrix would be the estimate obtained by the LR method.

\section{Conclusion}

We have proposed an extended flat-field model for spectral CT that exploits high correlation in the spectral flat-fields by replacing the spectral flat-fields with a low-rank approximation to mitigate ring artifacts. The proposed methodology can be combined with any existing reconstruction method and only depends on a single parameter, which is easy to choose by inspection of the singular values of the spectral flat-fields. The method was validated on neutron CT data and compared to the conventional flat-field correction and two preprocessing methods for ring reduction. Our method successfully mitigated ring artifacts in all experiments whereas the other methods struggled to suppress ring artifacts, especially in more challenging cases with severe ring artifacts.

\section*{Acknowledgments}{
This work was supported by The Villum Foundation (grant no.\ 25893). We gratefully acknowledge beamtime RB1820541 (DOI: 10.5286/ISIS.E.100529645) at the IMAT Beamline of the ISIS Neutron and Muon Source, Harwell, UK. EA was partially funded by the EPSRC grant EP/V007742/1 "Rich Nonlinear Tomography for Advanced Materials" and partially by the Federal Ministry of Education and Research (BMBF) and the Baden-Württemberg Ministry of Science as part of the Excellence Strategy of the German Federal and State Governments.}

\section*{References}
\bibliography{iopart-num}
\end{document}